\documentclass[12pt, twoside, a4paper]{article}
\usepackage[T1]{fontenc}
\usepackage[french]{babel}
\usepackage{amsmath,amsfonts,amssymb}
\newcommand{\prs}{\langle\;,\;\rangle}

\newcommand{\esp}{\quad\mbox{and}\quad}

\newcommand{\cohom}{\mathrm{H}}
\newcommand{\Om}{\Omega}
\newcommand{\om}{\omega}

\newcommand{\al}{\alpha}
\newcommand{\be}{\beta}

\newcommand{\Z}{\mathbb{Z}}

\newcommand{\N}{\mathbb{N}}

\newcommand{\R}{\mathbb{R}}

\newcommand{\Cinf}{\mathcal{C}^{\infty}}

\newcommand{\diff}{\mathrm{Diff}}

\newcommand{\too}{\longrightarrow}
\newcommand{\ga}{\gamma}
\newcommand{\Ga}{\Gamma}

\newcommand{\Lie}{\mathrm{Lie}}

\newcommand{\Isom}{\mathrm{Isom}}
\newcommand{\supp}{\mathrm{supp}}

\def \reel{ {\rm I}\!{\rm R} }

\def \Hyp{ { {\rm I}\!{\rm H}} }
\def \rat{ {\rm Q}\kern-.65em {}^{{}_/ }}

\def\ent{{{\rm Z}\mkern-5.5mu{\rm Z}}}

\newtheorem{theorem}{Theorem}[subsection]

\newtheorem{Corol}{Corollary}[subsection]

\newtheorem{remark}{Remark}[subsection]
\newtheorem{Lemma}{Lemma}[subsection]
\newtheorem{Example}{Example}[subsection]
\title{Cohomology of coinvariant differential forms}                                     \author{Abdelhak Abouqateb, Mohamed Boucetta and Mehdi Nabil } \date{ }\parindent=0cm
\begin{document}
\maketitle

\begin{abstract}
	Let $M$ be a smooth manifold and $\Gamma$ a group acting on $M$ by diffeomorphisms; which means that there is a  group morphism $\rho:\Gamma\rightarrow \diff(M)$ from $\Gamma$ to the group of diffeomorphisms  of $M$. For any such action we associate a cohomology $\cohom(\Omega(M)_\Gamma)$ which we
	call the cohomology of $\Gamma$-coinvariant forms. This is the cohomology of the graded vector space generated by the differentiable forms $\omega -\rho(\gamma)^*\omega$ where $\omega$ is a differential form with compact support and $\gamma\in \Gamma$. The present paper is an introduction to the study of this cohomology. More precisely, we study the relations between this cohomology, the de Rham cohomology and the cohomology of invariant forms $\cohom(\Omega(M)^\Gamma)$ in the case of isometric actions on compact Riemannian oriented manifolds and in the case of properly discontinuous actions on 
	manifolds.    %<------------------- 
	
\end{abstract}

\section{Introduction and statement of main results}\label{section1}

We start by fixing a notation which will be used along this paper. If $V$ is a real vector space,  $\Ga$ a group  and $\Gamma\times V\too V$, $(\ga,v)\mapsto \ga.v$ an action of  $\Gamma$ on $V$ by linear isomorphisms, we denote by $V^\Gamma$ the vector subspace of invariant vectors and $V_{\Gamma}$ the vector subspace of coinvariant vectors, i.e., $V_{\Gamma}=\mathrm{Span}\{v-\gamma.v,v\in V,\gamma\in\Gamma \}$.

Let $M$ be a connected differentiable manifold of dimension $n$ and $\Gamma$ a group. We denote by $\diff(M)$ the group of diffeomorphisms  of $M$. An action of $ \Gamma $ on $ M $ is  a group morphism $  \rho: \Gamma \rightarrow \diff (M) $.  The orbit space $ M / \Gamma $ will be endowed with its  quotient topology and we denote by  $\pi:M\longrightarrow M/\Gamma$ the canonical projection. For each $0\leq r\leq n$, $\Omega^r(M)$ is the space of $r$-differential forms on $M$ and $\Omega_c^r(M)$ the space of $r$-differential forms with compact support, where $\Omega_c^0(M)=\Cinf_c(M)$. For $\omega\in\Omega_c^r(M)$ and $\gamma\in \Gamma$, we denote by $\gamma^*\omega$ the pull-back of $\omega$ by the diffeomorphism $\rho(\gamma)$. The map $(\ga,\om)\too (\gamma^{-1})^*\omega$ defines an action of $\Ga$ both on 
$\Omega^r(M)$ and $\Omega_c^r(M)$. The graded vector spaces $\Omega^{}(M)^\Gamma:=\oplus_{r}\Omega^r(M)^\Gamma$ and $\Omega_c(M)_\Gamma:=\oplus_{r}\Omega_c^r(M)_\Gamma$ are stable under the usual de Rham differential operator $d$ and hence define two cohomologies $\cohom(\Omega(M)^\Gamma)$ and $\cohom(\Omega_c(M)_\Gamma)$. The cohomology $\cohom(\Omega(M)^\Gamma)$ known as the cohomology of invariant forms has been studied by many authors (for instance see \cite{G.H.V,  Haefliger, losik,tang}), however the cohomology $\cohom(\Omega_c(M)_\Gamma)$, we will call {\it {cohomology of $\Gamma$-coinvariant forms}}, is to our knowledge new and constitutes the main object of this paper.
The idea behind the introduction of this cohomology lies in the paper \cite{Abouqateb2} where the first author introduced the cohomology of $\mathcal{G}$-divergence forms and used it to study integrability conditions of an action of Lie algebra $\mathcal{G}$ on a manifold to a proper action of a connected Lie group. Our cohomology can  be seen, in a broad sense, as a sort of generalization of this idea in order to include actions of discrete groups or even actions of groups that are not necessarily Lie groups.

We shall denote also by $\cohom(M)$ (resp. $\cohom_c(M)$) the de Rham cohomology of $M$ (resp. the de Rham cohomology of differential forms with compact support). The group $\Ga$ acts linearly on $\cohom(M)$ and we call the elements of $\cohom(M)^\Ga$  invariant classes 	and the elements of $\cohom(M)_\Ga$ coinvariant classes. 

Through this paper, a cohomology class in $\cohom(M)$ will be denoted by $[\om]$, a cohomology class in $\cohom(\Omega(M)^\Gamma)$  by $[\om]^\Ga$ and a cohomology class in $\cohom(\Omega_c(M)_\Gamma)$  by $[\om]_\Ga$.

In this article, we  study in two situations the relations between the cohomologies $\cohom_c(M)$, $\cohom(\Omega(M)^\Gamma)$, $\cohom(\Omega_c(M)_\Gamma)$ and $\cohom(M)$. 
The first one is the case when  $\rho(\Gamma)\subset \Isom(M) $ where $\Isom(M)$ is the group of isometries of  a Riemannian metric on an oriented compact manifold $M$. Let $G^0$ be the connected component of the unity of $\Isom(M)$. Then $\rho(\Ga)\cap G^0$ is a normal subgroup of $\rho(\Ga)$ and the quotient $\Ga_1:=\rho(\Ga)/\rho(\Ga)\cap G^0$ is a finite group. Since the cohomology is invariant by homotopy, $\rho(\Ga)\cap G^0$ acts trivially on $\cohom(M)$ and hence the action of $\Ga$ on 
$\cohom(M)$ factors on an action of $\Ga_1$ on $\cohom(M)$. Our first main result is the following.
\begin{theorem}\label{main0} Let $\rho : \Gamma \rightarrow \Isom(M)$ be an action by isometries of a group $\Gamma$ on a compact oriented manifold. Then, for all $0\leq p \leq n$, the following assertions hold.
	\begin{enumerate}
		\item $\cohom^p(M)^{\Ga}=\cohom^p(M)^{\Ga_1}$, $\cohom^p(M)_{\Ga}=\cohom^p(M)_{\Ga_1}$ and $\cohom^p(M)=\cohom^p(M)^{\Ga_1}\oplus \cohom^p(M)_{\Ga_1}$,
		\item The map $\Phi:\cohom^p(\Omega(M)^\Gamma)\oplus \cohom^p(\Omega(M)_\Gamma)\too \cohom^p(M)$, $([\om]^\Ga+[\eta]_\Ga)\mapsto [\om+\eta]$ is an isomorphism, 
		\[ \Phi\left( \cohom^p(\Omega(M)^\Gamma)\right)=\cohom^p(M)^{\Ga_1}\esp
		\Phi\left( \cohom^p(\Omega(M)_\Gamma)\right)=\cohom^p(M)_{\Ga_1}.  \]
		\item If $\rho(\Ga)\subset G^0$ then $\cohom(M)$ is isomorphic to $ \cohom(\Omega(M)^\Gamma)$.
	\end{enumerate}
	
\end{theorem}

The last assertion of this theorem is an immediate consequence of the second one and, actually, has been proved by \'{E}lie Cartan  (see \cite{Haefliger}). One can also deduce from Theorem \ref{main0} that the cohomology of a compact connected Lie group $G$ is isomorphic to the cohomology of its Lie algebra by taking $M=G$, $\Ga=G$ and the action is by  left translations.

When a finite group $\Ga$ acts  on a manifold $M$ there exists always on $M$ an invariant Riemannian metric. So actions by finite groups on compact connected manifolds are particular cases of Theorem \ref{main0}.  The proof  in this case  is quite easy and does not involve the invariant metric and we will give it in the beginning of Section \ref{section2}. However, the proof in the general case needs more work and it will be based on
the following result which is an adaptation of the classical Hodge theory to the case of group actions. Note that classical Hodge theory has been used by many authors and for many purposes to study actions of Lie groups (see for instance \cite{alday, tang}).

\begin{theorem}[Hodge Decomposition Theorem of $\Omega(M)_\Gamma$ and $\Omega(M)^\Gamma$]  \label{main1} 
	Let $\rho : \Gamma \rightarrow \Isom(M)$ be an action by isometries of a group $\Gamma$ on a compact oriented manifold $M$. Then for all $0\leq p \leq n$, we  have the following orthogonal direct sum decompositions:
	\begin{enumerate}
		\item $\Omega^p(M)_\Gamma=d(\Omega^{p-1}(M)_\Gamma)\overset{\perp}{\oplus}\delta(\Omega^{p+1}(M)_\Gamma)\overset{\perp}{\oplus}\mathcal{H}^p(M)_\Gamma,$
		\item $\Omega^p(M)^\Gamma=d(\Omega^{p-1}(M)^\Gamma)\overset{\perp}{\oplus}\delta(\Omega^{p+1}(M)^\Gamma)\overset{\perp}{\oplus}\mathcal{H}^p(M)^\Gamma,$
		\item $\mathcal{H}^p(M)=\mathcal{H}^p(M)_\Gamma\overset{\perp}{\oplus}\mathcal{H}^p(M)^\Gamma$,
	\end{enumerate}
	where $\delta$ is the divergence operator and $\mathcal{H}^p(M)^\Gamma$ is the space of $\Gamma$-invariant harmonic $p$-forms and $\mathcal{H}^p(M)_\Gamma$ is the space of $\Gamma$-coinvariant harmonic $p$-forms.
\end{theorem}

The second case in our study is the case where the action of $\Gamma$ is  properly discontinuous and $M$ is not necessarily compact. 

A properly discontinuous action of a group $\Gamma$ on a manifold $M$ is an action such that,  for any compact $K\subset M$, the set $\Gamma_K=\{\gamma\in \Gamma\;,\;\gamma.K\cap K\neq\emptyset \}$ is finite. This  is equivalent to the map $ \Gamma{\times}M \ni (\gamma,x) \longmapsto (\gamma .x,x)\in M{\times}M$ being proper when we endow $\Gamma$ with the discrete topology. If $M$ is a connected Riemannian manifold, then every discrete subgroup of $\Isom(M)$ acts properly on $M$ (see \cite{KoT}). When $M$ is compact, an action $\rho : \Gamma \rightarrow \diff(M)$ is properly discontinuous if and only if the group $\Gamma$ is finite. For more details and reviews about properly discontinuous actions, one can see \cite{KoT, KoS}.

Suppose now that $\Gamma$ acts properly discontinuously on $M$.   We  define an {\it {average operator}} $$m:\Omega_c(M)\longrightarrow {\Omega(M)}^\Gamma \;\;\;\text{by}\;\; m(\omega)=\sum_{\gamma\in\Gamma}\gamma^*\omega.$$ It is easy to see that $m$ is a well-defined linear map because for any $x\in M$ the set $\Gamma_x=\{\gamma\in \Gamma\;,\;\gamma.x\in \rm{supp}\;\omega  \}$ is finite. Moreover, $m$  commutes with the differential $d$ and, for any $\omega\in\Omega_c(M)$, the form $m(\omega)$ is a $\Gamma$-invariant and has a $\Gamma$-compact support, which means that $\pi({\rm{supp}}\ m(\omega))=({\rm{supp}}\ m(\omega))/\Gamma$ is compact. We denote $\Omega(M)^\Gamma_{\Gamma c}$ the space of $\Gamma$-invariant forms on $M$ with $\Gamma$-compact support, this is a differential subalgebra of $\Omega(M)$. If $M/\Gamma$ is compact then $\Omega(M)^\Gamma_{\Gamma c}=\Omega(M)^\Gamma$.  Our second main result is the following theorem.

\begin{theorem}\label{main2}
	{{Let $\rho : \Gamma \rightarrow \diff(M)$ be a properly discontinuous action. Then  the operator $m:\Omega_c(M)\longrightarrow \Omega(M)^\Gamma_{\Gamma c}$ is surjective and $\ker{m}=\Omega_c (M)_\Gamma$. In other words, we obtain a short exact sequence of graded differential vector spaces:
			$$
			0\rightarrow \Omega_c (M)_\Gamma\overset{\iota}{\rightarrow}\Omega_c(M)\overset{m}{\rightarrow}\Omega(M)^\Gamma_{\Gamma c}\rightarrow 0,
			$$which gives rise to
			a long exact sequence in cohomology :
			$$
			\dots\rightarrow \cohom^p(\Omega_c (M)_\Gamma)\overset{\cohom(\iota)}{\longrightarrow} \cohom_c^p(M)\overset{\cohom(m)}{\longrightarrow} \cohom^p(\Omega(M)^\Gamma_{\Gamma c})\overset{\delta}{\longrightarrow}\cohom^{p+1}(\Omega_c (M)_\Gamma)\rightarrow\dots
			$$
			The connecting homomorphisms $\delta:\cohom^p(\Omega(M)^\Gamma_{\Gamma c})\longrightarrow \cohom^{p+1}(\Omega_c (M)_\Gamma)$ are given by the formula
			$$ \delta([\omega]^\Gamma)=[d\phi\wedge\omega]_\Gamma,$$
			where $\phi\in C^\infty(M)$ is a cutoff function. 
		}}\end{theorem}
		For the definition of  a cutoff function see Lemma \ref{coinvprop1}. Moreover, a cutoff function satisfies $\sum_{\ga\in\Ga} \phi\circ\rho(\ga)=1$ and hence  $m(d\phi\wedge\om)=0$ which shows that $d\phi\wedge\om\in \Om_c^{p+1}(M)_\Ga$ and the expression of $\delta$ given in the theorem makes sense.
		
		The paper is organized as follows. In Section \ref{section2}, we prove Theorems \ref{main0} and \ref{main1}. We devote Section \ref{section3} to a proof of Theorem \ref{main2} and some of its applications.

		\section{Actions by isometries on compact manifolds}\label{section2}
		
		This section is devoted to proving Theorem \ref{main0} and \ref{main1}. The proof of Theorem \ref{main0} is based on Theorem \ref{main1}. However, in the case when $\Ga$ is finite we can prove Theorem \ref{main0} straight forwardly. This proof presents some interest and we will give it first. 
		
		Let us start with a general remark. Let $V$ be a real vector space and $\Gamma\times V\too V$, $(\ga,v)\mapsto \ga.v$ an action of a finite group $\Gamma$ on $V$ by linear isomorphisms. The {\it average operator} is the linear map $m:V\too V^\Ga$ given by
		\begin{equation}\label{m}
		m(v)=\frac{1}{|\Gamma|}\sum_{\gamma\in\Gamma}\gamma.v,
		\end{equation}where $|\Gamma|$ is the cardinal of $\Ga$. It is an easy task to check that
		\begin{equation}\label{eq1}
		v -m(v)=\sum_{\gamma\in\Gamma} (\dfrac{v}{|\Gamma|}-\gamma. (\dfrac{v}{|\Gamma|})).
		\end{equation}This formula shows that
		\begin{equation}\label{eq2}
		\ker m=V_\Ga\esp V=V^\Ga\oplus V_\Ga.
		\end{equation}
		
		A non trivial and analogous splitting to \eqref{eq2} has been obtained in a more general setting involving
		linear actions of compact Lie groups on some classes of topological vector spaces (see \cite[Theorem 3.36 pp. 76]{hof}).

		\begin{remark}\label{rem} Note that if $\Ga$ is a group (not necessarily finite) acting by isometries on a real pre-Hilbert space $(V,\langle\;,\;\rangle)$ then $V^\Ga=(V_\Ga)^\perp$ and hence if $\dim V$ is finite then $V=V^\Ga\oplus V_\Ga$.
			
		\end{remark}
		
		{\bf Proof of Theorem \ref{main0} in case where $\Ga$ is finite.} 
			According to \eqref{eq2}, for any $r\in\{0,\ldots,n\}$,
			\begin{equation} \label{eq40}
			\Omega^r(M)=\Omega^r(M)^\Gamma\oplus \Omega^{r}(M)_\Gamma\esp \cohom^r(M)=\cohom^r(M)^\Ga\oplus \cohom^r(M)_\Ga.
			\end{equation}
			An immediate consequence of the first relation is that $\Phi$ is an isomorphism. Moreover, it is obvious that $\Phi(\cohom^r(\Omega(M)^\Gamma))\subset  \cohom^r(M)^\Ga$. Conversely, let $[\om]\in \cohom^r(M)^\Ga$. Then $\om$ is closed and, for any $\ga\in\Ga$, there exists $\al^\ga$ such that $\ga^*\om=\om+d\al^\ga$. It follows then that $m(\om)=\om+d\left( \frac1{|\Ga|}\sum\al^\ga\right)$ and hence $[\om]=\Phi([m(\om)]^\Ga)$. Thus $\Phi(\cohom^r(\Omega(M)^\Gamma))=  \cohom^r(M)^\Ga$. On the other hand, $\cohom^r(M)_\Ga\subset \Phi(\cohom^r(\Omega(M)_\Gamma))$. Indeed, if $[\om]\in \cohom^r(M)_\Ga$, then there exists a family $(\om_i)$ of closed forms  and a family $(\ga_i)$ of elements of $\Ga$ and a form $\al$ such that
			\[ \om=\sum_{i}\left(\om_i-\ga_i^*\om_i\right)+d\al. \]
			This shows that $[\om]=\Phi([\sum_{i}\left(\om_i-\ga^*\om_i\right)]_\Ga)$. Now since $M$ is compact, its cohomology is finite dimensional. Then from the second relation in \eqref{eq40} and the fact that $\Phi(\cohom^r(\Omega(M)^\Gamma))=  \cohom^r(M)^\Ga$, we deduce that
			$\dim \cohom^r(M)_\Ga =\dim\Phi(\cohom^r(\Omega(M)_\Gamma))$ which completes the proof of Theorem \ref{main0} in the case when $\Ga$ is finite.
		$\square$
		
		We will prepare now the proof of Theorem \ref{main1} by recalling the main results of Hodge theory of harmonic forms on a Riemannian manifold $M$. For more details  one can see \cite{Warner}.
		
		We suppose now that $M$ is a connected compact oriented Riemannian manifold of dimension $n$. 
		We denote by $\prs$ the Riemannian metric and 
		by $\mu$
		the Riemannian volume form associated to it.  For every $p\in \{0,\ldots,n\}$, let $\Lambda^pT^*M$ be the vector bundle of alternating covariant $p$-tensors on $M$. The Riemannian metric induces a Riemannian metric on the vector bundle $\Lambda^pT^*M\rightarrow M$ denoted  by the same notation. 
		For any $p\in\N$, we consider the  $\ast$-operator which is the bundle isomorphism
		$\ast : \Lambda^pT^*M  \longrightarrow \Lambda^{n-p}T^*M$  uniquely determined by the relation
		\[ \al\wedge *\be=\langle\al,\be\rangle\mu,\; \al,\be\in \Lambda^pT^*M. \]		 
		The star operator extends to a mapping $ \ast : \Omega^p(M) \rightarrow \Omega^{n-p}(M)$   called the {\it {Hodge star  operator}}. The {\it {codifferential operator}} on $p$-forms $\delta:\Omega^p(M)\longrightarrow \Omega^{p-1}(M)$ is defined by $\delta\omega=(-1)^{n(p+1)+1}*d*\omega$ and the {\it {Laplace-Beltrami operator }}is  defined as $\Delta:\Omega^p(M)\longrightarrow\Omega^p(M)$ by $ \Delta=\delta d+d\delta$. A differential form $\omega$ is said to be  harmonic if $\Delta \omega = 0$. We denote by $\mathcal{H}^p(M)$ the space of $p$-harmonic forms with $\mathcal{H}(M):=\bigoplus\limits_{p=0}^{n}\mathcal{H}^p(M)$. The last major definition is the scalar product $(\cdot,\cdot )_M: \Omega^p(M)\times\Omega^p(M)\longrightarrow \R$ given by 
		$$( \alpha, \beta )_M=\int_M\alpha\wedge \ast\beta. 
		$$
		We extend $( \ . , . )_M$ to an inner product on $\Omega^\ast(M)$ by declaring $\Omega^p(M)$ and $\Omega^q(M)$ to be orthogonal for $p\neq q$. We denote the corresponding norm by $|| \omega ||$. There are the main properties of these operators:
		\begin{enumerate}
			\item[$(i)$] The Laplace-Beltrami $\Delta$ and Hodge star operator commute : $*\Delta=\Delta*$.
			\item[$(ii)$] The codifferential is the adjoint of the exterior derivative with respect to $(\cdot,\cdot )_M$, i.e., for all $\alpha,\ \beta$ we have $$( d\alpha,\beta)_M=( \alpha,\delta\beta)_M.$$
			\item[$(iii)$]  The Laplace-Beltrami operator is self-adjoint $( \Delta\alpha,\beta)_M=( \alpha,\Delta\beta)_M$.
			\item[$(iv)$] For every $\alpha\in\Omega(M)$, $\Delta\alpha=0$ if and only if $d\alpha=0$ and $\delta\alpha=0$.
		\end{enumerate}	
		The main result in this theory is the following theorem.
		\begin{theorem}[Hodge Decomposition Theorem] 
			For all $0\leq p \leq n$, $\mathcal{H}^p(M)$ is finite dimensional and we have the following orthogonal direct sum decomposition of $\Omega^p(M)$:
			\begin{eqnarray}\label{h}
			\Omega^p(M)
			&=& d(\Omega^{p-1}(M))\overset{\perp}{\oplus}\delta(\Omega^{p+1}(M))\overset{\perp}{\oplus}\mathcal{H}^p(M).
			\end{eqnarray} Moreover, every de Rham cohomology class on a compact oriented Riemannian manifold $M$ has a unique harmonic representative. More precisely, the mapping $\jmath:\mathcal{H}^p(M)\longrightarrow \cohom^p(M)$ given by: $\jmath(\omega)=[\omega]$ is an isomorphism of vector spaces.
		\end{theorem}
		
		\paragraph{\bf Proof of Theorem \ref{main1}}
		
			Since $\rho(\Ga)\subset \Isom(M)$, for any $\ga\in\Ga$, $\ga^*$ commutes with $\delta$  and $\Delta$, $\ga^*$ is an isometry of $\Om(M)$ endowed with the scalar product $(.,.)_M$. Thus, for any $p\in\{0,\ldots,n\}$, $\Ga$ acts by isometry on $\mathcal{H}^p(M)$. Since $\dim \mathcal{H}^p(M)_\Ga$ is finite and according to Remark \ref{rem}, we have
			\begin{equation}\label{eq12}
			\mathcal{H}^p(M)=\mathcal{H}^p(M)^\Ga\oplus \mathcal{H}^p(M)_\Ga.
			\end{equation} This shows the third relation in the theorem. 
			
			Before proving the first and the second decomposition, note that  $\mathcal{H}^p(M)_\Ga=\mathcal{H}^p(M)\cap \Om^p(M)_\Ga$. Indeed, if $\al-\ga^*\al\in \mathcal{H}^p(M)$ then,  by virtue of \eqref{h},   $\al=d\al_1+\delta\al_2+\al_3$ where $\al_3$ is harmonic. Hence $\al-\ga^*\al=d(\al_1-\ga^*\al_1)+\delta(\al_2-\ga^*\al_2)+(\al_3-\ga^*\al_3)$. Since $\al_3-\ga^*\al_3$ is harmonic, we deduce that 
			$\al-\ga^*\al=\al_3-\ga^*\al_3$ and the result follows.

			The first decomposition is quite easy. Indeed,	let $\omega\in\Omega^p(M)$ and $\gamma\in\Gamma$. By virtue of \eqref{h},
			$\omega=d\alpha+\delta\beta+\eta,$ where  $\eta\in\mathcal{H}^p(M)$. Then 
			$$ \omega-\gamma^*\omega=d(\alpha-\gamma^*\alpha)+\delta(\beta-\gamma^*\beta)+(\eta-\gamma^*\eta). $$
			Since $\eta-\gamma^*\eta$ is harmonic we get the first decomposition.
			
			The second decomposition needs some work. Let $\omega\in\Omega^p(M)^\Gamma$ and as  previously we write
			\begin{equation}
			\label{hareq1}
			\omega=d\alpha+\delta\beta+\eta,
			\end{equation}
			where $\alpha\in\Omega^{p-1}(M)$, $\beta\in \Omega^{p+1}(M)$ and $\eta\in\mathcal{H}^p(M)$. Since $\omega$ is $\Gamma$-invariant we obtain that for any $\gamma\in\Gamma$, 
			$$ \omega-\gamma^*\omega=0=d(\alpha-\gamma^*\alpha)+\delta(\beta-\gamma^*\beta)+\eta-\gamma^*\eta.$$
			which implies that $d\alpha=\gamma^*d\alpha$, $\delta\beta=\gamma^*\delta\beta$ and $\eta=\gamma^*\eta$ for any $\gamma\in\Gamma$. In particular, $\eta\in\mathcal{H}^p(M)^\Gamma$. We shall now show that we can replace $\alpha$ and $\beta$ by $\Gamma$-invariant differential forms, more precisely we will show that
			$\alpha=\alpha_1+\alpha_2$ and $\beta=\beta_1+\beta_2,$
			where $\alpha_1\in\Omega^{p-1}(M)^\Gamma$, $\beta_1\in \Omega^{p+1}(M)^\Gamma$, $d\alpha_2=0$ and $\delta\beta_2=0$. To do so, we write $\alpha=d\mu+\delta\nu+\lambda$ and  $\beta=d\hat{\mu}+\delta\hat{\nu}+\hat{\lambda}$ the respective Hodge decompositions of $\alpha$ and $\beta$. Hence, for any $\gamma\in\Gamma$,
			$$ d(\alpha-\gamma^*\alpha)=d\delta(\nu-\gamma^*\nu)\esp \delta(\beta-\gamma^*\beta)=\delta d(\hat{\mu}-\gamma^*\hat{\mu}).$$
			This implies that $d\delta(\nu-\gamma^*\nu)=0$ and $\delta d(\hat{\mu}-\gamma^*\hat{\mu})=0$. Therefore
			$$( \delta(\nu-\gamma^*\nu),\delta(\nu-\gamma^*\nu) )_M=( \nu-\gamma^*\nu,d\delta(\nu-\gamma^*\nu) )_M=0.$$ 
			Thus $\delta\nu-\gamma^*\delta\nu=0$ for any $\gamma\in\Gamma$ and hence $\delta\nu\in\Omega^{p-1}(M)^\Gamma$. In the same way, we can show that $d\hat{\mu}\in\Omega^{p+1}(M)^\Gamma$. We put $\alpha_1=\delta\nu$, $\alpha_2=d\mu+\lambda$, $\beta_1=d\hat{\mu}$ and $\beta_2=\delta\hat{\nu}+\hat{\lambda}$. By replacing in \eqref{hareq1}, we obtain
			$$ \omega=d(\alpha_1+\alpha_2)+\delta(\beta_1+\beta_2)+\eta=d\alpha_1+\delta\beta_1+\eta, $$
			and we have $\alpha_1\in\Omega^{p-1}(M)^\Gamma$, $\beta_1\in \Omega^{p+1}(M)^\Gamma$ and $\eta\in\mathcal{H}^p(M)^\Gamma$.
			This completes the proof.
		$\square$

		{\bf Proof of Theorem \ref{main0}}. 
			The facts that $\cohom^p(M)_{\Ga}=\cohom^p(M)_{\Ga_1}$ and $\cohom^p(M)^{\Ga}=\cohom^p(M)^{\Ga_1}$ are obvious and the decomposition $\cohom^p(M)=\cohom^p(M)_{\Ga_1}\oplus \cohom^p(M)^{\Ga_1}$ is a consequence of the fact that $\Ga_1$ is finite and \eqref{eq2}. From the first and the second decomposition in Theorem \ref{main1}, we get that there is an isomorphism
			\[ \Phi_1:\cohom^p(\Omega(M)^\Gamma)\oplus \cohom^p(\Omega(M)_\Gamma)\too \mathcal{H}^p(M)^\Ga\oplus \mathcal{H}^p(M)_\Ga\stackrel{\eqref{eq12}}{=}\mathcal{H}^p(M), \]and $\Phi=j\circ\Phi_1$ where $j:\mathcal{H}^p(M)\too \cohom^p(M)$ is the natural isomorphism. Since $j(\mathcal{H}^p(M)^\Ga)=\cohom^p(M)^\Ga$ and $j(\mathcal{H}^p(M)_\Ga)=\cohom^p(M)_\Ga$, we get that
			\[ \Phi\left( \cohom^p(\Omega(M)^\Gamma)\right)=\cohom^p(M)^{\Ga_1}\esp
			\Phi\left( \cohom^p(\Omega_c(M)_\Gamma)\right)=\cohom^p(M)_{\Ga_1}.  \]
			If $\rho(\Ga)\subset G^0$ then $\Ga_1$ is trivial and we get that the cohomology of $M$ is isomorphic to the cohomology of $\Ga$-invariant forms.$\square$

		\begin{remark} Under the hypothesis of Theorem \ref{main0}, we have from Remark \ref{rem} that for any $p\in\{0,\ldots,n\}$, $(\Om^p(M)_\Ga)^\perp=\Om^p(M)^\Ga$ and hence
			\[ (\Om^p(M)^\Ga\oplus \Om^p(M)_\Ga)^\perp=\{0\}. \]However, we can not derive any useful conclusion which could permit us to avoid the use of Hodge theory in the proof of the theorem.
			
		\end{remark}

		\section{Cohomology of properly discontinuous actions}\label{section3}
		
		In this section, we prove Theorem \ref{main2}. We refer the reader to the introduction where the definitions of properly discontinuous actions 	and the average operator associated to these actions have been given. 
		
		The main tool we will need is the following lemma about the existence of cutoff functions for properly discontinuous actions. These functions were already discussed in Bourbaki \cite[Proposition 8, p. 51]{Bourbaki}  in the context of topological group actions, and in the article of El Kacimi and Matsumoto \cite{KM} in the case when the projection $M\rightarrow M/\Gamma$ is a covering, our proof however does not presuppose that the action is free and is therefore available in a more general setting. 
		\begin{Lemma}[Cutoff functions]
			\label{coinvprop1}
			Let $\Gamma$ be a discrete group acting properly discontinuously on a manifold $M$. Then there exists a $C^\infty$ positive function $\phi:M\longrightarrow \R$  such that for any compact $B\subset M/\Gamma$, $\supp(\phi)\cap\pi^{-1}(B)$ is compact. Furthermore, we have 
			\begin{equation}
			\label{eq0}
			\sum\limits_{\gamma\in\Gamma}\phi\circ\gamma=1.
			\end{equation}
		\end{Lemma}
		{\bf Proof.}
			Let $(V_n)_{n\in\N}$ be a locally finite covering of $M/\Gamma$ by relatively compact open sets. We claim that there exists a locally finite covering $(W_n)_n$ of $M/\Gamma$ such that $W_n$ is relatively compact and $\overline{V}_n\subset W_n$. Indeed, denote by $J_1=\{i\in\N,\;\overline{V}_1\cap V_i\neq\emptyset\}$. Since $\overline{V}_1$ is compact and $(V_n)_n$ is a locally finite cover of $M/\Gamma$, we obtain that $J_1$ is finite and moreover 
			$$\overline{V}_1\subset\underset{j\in J_1}{\bigcup} V_j\;\;\;\text{and}\;\;\;\overline{V}_1\cap\underset{j\notin J_1}{\bigcup} V_j=\emptyset.$$
			There exists then a relatively compact open set $W_1$ such that $\overline{V}_1\subset W_1$ and 
			$$W_1\subset\underset{j\in J_1}{\bigcup} V_j\;\;\;\text{and}\;\;\;W_1\cap\underset{j\notin J_1}{\bigcup} V_j=\emptyset.$$
			Hence, the family $\mathcal{F}_1=\{W_1\}\cup\{V_j,\;j\geq 2\}$ is a locally finite cover of $M/\Gamma$ by relatively compact open sets. Repeating this process on the family $\mathcal{F}_1$, we prove by induction our claim.
			
			Furthermore, for all $n\in \N$, there exists two relatively compact open sets $U_n$ and $O_n$ of $M$ satisfying $\pi(U_n)=V_n$, $\pi(O_n)=W_n$ and $\overline{U}_n\subset O_n$. Indeed, we start by taking a relatively compact open set $\hat{O}_n$ such that $\pi(\hat{O}_n)=W_n$  and $\hat{U}_n$ a relatively compact open set of $M$ such that $\pi(\hat{U}_n)=V_n$. Since $V_n\subset W_n$ then $\hat{U}_n\subset\bigcup_{\gamma\in\Gamma}\gamma \hat{O}_n$. We then put :
			$$ U_n=\hat{O}_n\cap\underset{\gamma\in\Gamma}{\bigcup}\gamma^{-1}\hat{U}_n.$$
			Thus $U_n$ is a relatively compact open set of $M$ and it is clear that $\pi(U_n)\subset V_n$. Conversely, if $x\in V_n$, then there exists $a\in\hat{U}_n$ such that $\pi(a)=x$. Since $\hat{U}_n\subset\bigcup_{\gamma\in\Gamma}\gamma \hat{O}_n$, there exists $\gamma\in\Gamma$ such that $a\in\gamma \hat{O}_n$ and hence $\gamma^{-1}a\in\gamma^{-1}\hat{U}_n\cap \hat{O}_n$, moreover $\pi(\gamma^{-1}a)=x$. Thus, we  obtain that $x\in\pi(U_n)$ which means that $\pi(U_n)=V_n$. Now denote $\{\gamma_{1}^n,\dots,\gamma_{r_n}^n\}=\{\gamma\in\Gamma,\;\overline{U}_n\cap\gamma\hat{O}_n\neq\emptyset\}$. Since $\overline{V}_n\subset W_n$, then $\overline{U}_n\subset\bigcup_{\gamma\in\Gamma}\gamma\hat{O}_n$. We deduce that 
			$$ \overline{U}_n\subset \overset{r_n}{\underset{i=1}{\bigcup}}\gamma\hat{O_n}:=O_n.$$
			To summarize, $O_n$is a relatively compact open set of $M$ satisfying $\pi(O_n)=W_n$ and $\overline{U}_n\subset O_n$.
			
			After this tedious construction, we can begin the proof of the lemma. Define for all $n\in\N$ the function $g_n\in C^\infty_c(M)$ satisfying 
			\begin{enumerate}
				\item $0\leq g_n\leq 1$ et $\supp(g_n)\subset O_n$,
				\item $g_n=1$ on $U_n$.
			\end{enumerate}
			Put $g=\sum_{n\in\N} g_n$. Then $g$ is well-defined, positive and of class $C^\infty$ on $M$. Indeed, let $U$ be a relatively compact open set on $M$ and set $J=\{n\in \N,\overline{U}\cap O_n\neq\emptyset\}.$ Since $(W_n)_{n\in\N}$ is a locally finite cover of $M/\Gamma$, then  $(O_n)_{n\in\N}$ is a locally finite cover of $M$ and thus the set $J$ is finite. Hence 
			$$ g_{_{|U}}=\sum_{j\in J} {g_j}_{_{|U}}\in C^\infty(U).$$
			This shows that $g$ is well defined and of class $C^\infty$. Let us next consider a compact set $B$ de $M/\Gamma$ and put $I_0=\{n\in \N,\;\pi^{-1}(B)\cap O_n\neq \emptyset\}$. Since $(W_n)_n$ is a locally finite cover of $M/\Gamma$ and $\pi(O_n)=W_n$, we have 
			$$ I_0\subset\{n\in \N,\;\pi^{-1}(B)\cap \pi^{-1}(W_n)\neq \emptyset\}=\{n\in \N,\;B\cap W_n\neq \emptyset\}\;\;\text{which is finite}.$$
			Thus $\pi^{-1}(B)\cap\supp(g)\subset\pi^{-1}(B)\cap\bigcup_{n\in\N}O_n=\pi^{-1}(B)\cap\bigcup_{i\in I_0} O_i$. Since $\bigcup_{i\in I_0} O_i$ is relatively compact, we obtain that $\pi^{-1}(B)\cap\supp(g)$ is compact for any compact $B$ in $M/\Gamma$. Next let $x\in M$, there exists then $n\in\N$ such that $\pi(x)\in V_n$. This means that there is a $\gamma\in\Gamma$ such that $\gamma x\in U_n$ and thus $g(\gamma x)>0$. Hence $\sum_{\gamma\in\Gamma}g(\gamma x)>0$ for all $x\in M$. Finally, we put  
			$$\phi=\frac{g}{\sum_{\gamma\in\Gamma} g\circ\gamma}.$$
			This gives that $\supp(\phi)=\supp(g)$ and $\sum_{\gamma\in\Gamma}\phi\circ\gamma=1$. This ends the proof.$\square$
		
		\begin{remark}
			When $M/\Gamma$ is compact, the cutoff function $\phi$ in Lemma \ref{coinvprop1} has a compact support.
		\end{remark} 
		{\bf Proof of Theorem  \ref{main2}. }
		
		 We start by showing that $m$ is onto and $\ker m=\Omega_c(M)_\Gamma$.
			Let $\eta\in\Omega(M)^\Gamma_{\Gamma c}$ and choose a cutoff function $\phi\in C^\infty(M)$ as in Lemma \ref{coinvprop1}. The set $K=\pi(\supp\;\eta)$ is compact and since $\supp(\eta)$ is $\Gamma$-invariant then $\supp(\eta)=\pi^{-1}(K)$. It follows from the properties of $\phi$ that  $\omega=\phi\eta$ has compact support in $M$ since $\supp\;\omega\subset\supp\;\phi\cap\pi^{-1}(K)$ which is compact. Moreover,
			$$ m(\omega)=\sum_{\gamma\in\Gamma}\gamma^*\omega=\sum_{\gamma\in\Gamma}\gamma^*(\phi\eta)=\left(\sum_{\gamma\in\Gamma}\gamma^*\phi\right)\eta=\eta.$$
			This shows that $m$ is onto. On the other hand, for any $\alpha\in\Gamma$ and any $\omega\in\Omega_c(M)$, we have 
			$$ m(\omega-\alpha^*\omega)=\sum_{\gamma\in\Gamma}\gamma^*\omega-\sum_{\gamma\in\Gamma}\gamma^*\alpha^*\omega=\sum_{\gamma\in\Gamma}\gamma^*\omega-\sum_{\gamma\in\Gamma}\gamma^*\omega=0.$$
			Hence $\Omega_c(M)_\Gamma\subset\ker(m)$. Conversely, let $\omega\in\ker(m)$. Put $K_0=\supp\;\om$ and for all $\gamma\in\Gamma$ 
			$$\phi_\gamma=\phi\circ\rho(\gamma)\esp \omega_\gamma=\phi_\gamma\omega-(\gamma^{-1})^*(\phi_\gamma\omega)=\phi_\gamma\omega-
			\phi(\gamma^{-1})^*\omega.$$
			It is clear that, for all $\gamma\in\Gamma$,  $\omega_\gamma\in\Omega_c^*(M)_\Gamma$. Put
			$\Ga_\om=\left\{ \gamma\in\Gamma:\om_\ga\not=0    \right\}$. We claim that $\Ga_\om$ is finite.
			Indeed,  $\pi(K_0)$   is a compact set contained in $M/\Gamma$ and then $K_1=\supp\;\phi\cap\pi^{-1}(\pi(K_0))=\bigcup_{\gamma\in\Gamma}(\supp\;\phi\cap\gamma K_0)$ is a compact set in $M$. Since the action of $\Ga$ is properly discontinuous and
			\[ A=\{\gamma\in\Gamma,\;K_1\cap(\gamma K_0)\neq\emptyset\}
			\subset \{\gamma\in\Gamma,\;(K_1\cup K_0)\cap(\gamma (K_1\cup K_0))\neq\emptyset\}, \]it follows that $A$ is finite. Moreover, $A=\{\gamma\in\Gamma,\;\supp\;\phi\cap\gamma K_0\neq\emptyset\}$.
			Hence, if $\ga\in\Gamma\setminus A$, it follows that  $$\supp(\phi_\ga\omega)\subset(\ga^{-1}\supp\;\phi)\cap\supp\;\omega=\ga^{-1}(\supp\;\phi\cap\ga\supp\;\omega)=\emptyset.$$
			Consequently $\phi_\ga\omega=0$ and thus $\omega_\ga=0$ for all $\ga\in\Gamma\setminus A$. This  shows the claim. It follows that
			\begin{eqnarray*}
				\sum_{\gamma\in \Ga_\om}\omega_\gamma&=&\sum_{\gamma\in \Gamma}\omega_\gamma=\sum_{\gamma\in \Gamma}\phi_\gamma\omega-\sum_{\gamma\in \Gamma}\phi(\gamma^{-1})^*\omega\\&=&\left(\sum_{\gamma\in \Gamma}\phi_\gamma\right)\omega-\phi\left(\sum_{\gamma\in \Gamma}(\gamma^{-1})^*\omega\right)\\&=&\omega-\phi m(\omega)=\omega.\end{eqnarray*}
			Thus $\omega=\sum_{\gamma\in A}\omega_\gamma\in\Omega_c(M)_\Gamma$ and completes of the proof of $\ker(m)=\Omega_c(M)_\Gamma$.
			Hence, we get  a short exact sequence $$
			0\rightarrow \Omega_c(M)_\Gamma\overset{\iota}{\rightarrow}\Omega_c(M)\overset{m}{\rightarrow}
			\Omega(M)^\Gamma_{\Gamma c}\rightarrow 0.
			$$ The existence of a long exact sequence in cohomology is a consequence of a well-known result.
			
			Now, let $\omega\in\Omega^p(M)^\Gamma$ be a closed form. 
			From the expression of the connecting homomorphism in the snake lemma we have $\delta([\omega])=[\beta]_\Gamma$ where $\beta=d\alpha$ and $\alpha\in \Omega_c^p(M)$ is such that $m(\alpha)=\omega$. We choose then $\alpha=\phi\omega$, which gives $\beta=d\phi\wedge\omega$. This completes the proof of the theorem.$\square$

		Before giving some corollaries to Theorem \ref{main2} and some examples to illustrate it, let's start by some remarks on $\cohom^{0}(\Omega_c (M)_\Gamma)$ and $\cohom^{1}(\Omega_c (M)_\Gamma)$ in the case of properly discontinuous actions.

		Let $\rho : \Gamma \rightarrow \diff(M)$ be a properly discontinuous action.
		\begin{enumerate}
			\item Note first that $\cohom^{0}(\Omega_c (M)_\Gamma)=0$. Indeed, if $[f]_\Gamma\in \cohom^{0}(\Omega_c^\ast (M)_\Gamma)$, then $f$ is a constant and $\sum_{\gamma\in\Gamma}\gamma^*f=0$. If $M$ is not compact, then $f=0$ because it is constant function with compact support. If $M$ is compact, then $\Gamma$ is finite group and hence $\sum_{\gamma\in\Gamma}\gamma^*f=|\Gamma|f=0$, that is $f=0$.
			\item Suppose now that $M/\Ga$ is compact. Then for every cutoff function $\phi$, we have that $\phi$ has a compact support and $m(\phi)=1$. Thus $m(d\phi)=0$ and hence $d\phi\in \Omega_c^1 (M)_\Gamma$. So $[d\phi]_\Ga$ defines a cohomology class in $\cohom^{1}(\Omega_c (M)_\Gamma)$ and one can see easily that this class does not depend on the choice of $\phi$. We denote this class by $\theta_\Ga$. Obviously if $\cohom(\iota):\cohom^{1}(\Omega_c (M)_\Gamma)\rightarrow \cohom_c^1(M) $ is the natural map, then $\cohom(\iota)(\theta_\Ga)=0$. 
			Moreover, if $M$ is non-compact, then $\theta_\Ga\not=0$. Indeed, if there exists a coinvariant function $\xi$ with compact support and $\phi-\xi$ is constant, one must have $\phi-\xi=0$
			since $\phi-\xi$ is a constant function with compact support in a non-compact manifold. Thus $\cohom^{1}(\Omega_c (M)_\Gamma)\not=\{0\}$ and $\cohom(\iota)$ is not injective.
			\item Suppose now that $M/\Ga$ is non-compact. In this case $\cohom(\iota)$ is always injective. Indeed, let $\omega\in\Omega_c(M)_\Gamma$ be a closed form such that $[\omega]=0$. Write $\omega=df$ such that $f\in C^\infty_c(M)$. Since we have $dm(f)=m(df)=m(\omega)=0$ we obtain that $m(f)$ is a constant function with $\Gamma$-compact support. If $m(f)\neq 0$ then $\supp\;m(f)=M$ and consequently $\supp\;m(f)/\Gamma=M/\Gamma$ which is noncompact, this leads to a contradiction. Thus $m(f)=0$ which is equivalent to $f\in C^\infty_c(M)_\Gamma$, hence $[\omega]_\Gamma=0$. We conclude that $\iota$ is injective.
		\end{enumerate}

		Let us give now an important consequence of the long exact sequence in Theorem \ref{main2}.
		
		\begin{Corol}\label{co}
			Let $M$ be a contractible manifold and $\rho : \Gamma \rightarrow \diff(M)$ be a properly discontinuous action with compact orbit space $M/\Gamma$. Then for any $1\leq p\leq  n$, we have  $$\cohom^p(\Omega_c(M)_\Gamma)\simeq \cohom^{p-1}(\Omega(M)^\Gamma).$$ In particular, $\cohom^1(\Omega_c(M)_\Gamma)=\rm{Span\{[\theta]_\Gamma}\}$.
		\end{Corol}
		
		We end this paper by giving some examples which illustrates Theorem \ref{main2}. These examples show that, in contrast with the case studied in Theorem \ref{main0}, for properly discontinuous actions the cohomology $\cohom_c(M)$ could be trivial but the cohomologies $\cohom(\Omega(M)^\Gamma)$ and $\cohom(\Omega(M)_\Gamma)$ can be interesting.
		
		\begin{Example}\begin{enumerate}

				\item 	{\bf Nilmanifolds.} A compact nilmanifold is the quotient of a simply connected nilpotent Lie group $G$ by a discrete subgroup $\Gamma$ of $G$ such that $G/\Gamma$ is compact $($\cite{Ragunathan}$)$. The simplest example being the $n$-dimensional torus viewed as a quotient of $\reel^n$ by $\ent^n$. The cohomology of $\Gamma$-invariant forms on $G$ can naturally be identified with the cohomology of $G/\Gamma$ and a famous theorem of Nomizu (\cite{N}) asserts that the cohomology of $\cohom(G/\Gamma)$ is isomorphic to $\cohom(\mathfrak{g})$ the cohomology of the Lie algebra $\mathfrak{g}=\Lie(G)$. So we get from Corollary \ref{co} that, for every $p\in\{1,\ldots,n\}$,
				\[ \cohom^p(\Omega(G)^\Gamma)\simeq \cohom^{p}(\mathfrak{g})\esp \cohom^p(\Omega_c(G)_\Gamma)\simeq \cohom^{p-1}(\mathfrak{g}). \]	 
				In the particular case of the usual action by translations of $\Z^n$ on $\R^n$, we obtain $\dim \cohom^p(\Omega_c(\reel^n)_{\Z^n})=C^{p-1}_n$.
				
				\item 	{\bf Riemann Surfaces.} 
				Let $\Sigma_g$ be a connected compact Riemann surface of genus $g\geq 2$. The fundamental group of $\Sigma_g$ can be identified with a discrete subgroup $\Gamma_g$ of $\rm{PSL}(2,\reel)=\rm{SL}(2,\reel)/{\{\pm\rm{I}\}}$ so that the surface $\Sigma_g$ is identified with the orbit space $\Hyp/{\Gamma_g}$ of the action of $\Gamma_g$ on the Poincar\'e half-plane $\Hyp$. This action is given by
				$$(A,z)\longmapsto \frac{az+b}{cz+d}
				$$
				for every $z\in \Hyp$ and $
				\left(\begin{array}{cc} a & b \\ c & d
				\end{array}\right)$ a matrix of $\rm{SL}(2,\reel)$ representative $A\in\Gamma_g$. The $2$-cohomology space of $\Gamma_g$-coinvariant forms on $\Hyp$ is then isomorphic to the $1$-cohomology space of $\Sigma_g$, this gives
				$\dim(\cohom^2(\Omega_c(\Hyp)_{\Gamma_g}) = 2g.
				$
				\item {\bf Compact Clifford-Klein forms.} 
				A homogeneous space $G/H$ is said to have a compact Clifford-Klein form if there exists a discrete subgroup $\Gamma$ of $G$ which acts properly discontinuously on $G/H$, such that the quotient space $\Gamma\backslash G/ H$ is compact. The double coset space $\Gamma\backslash G/ H$ is then called a compact Clifford-Klein form. Compact Clifford-Klein forms has been studied by many authors (see\cite{Ko} for instance).
				
				For every compact Clifford-Klein form $\Gamma\backslash G/ H$, by virtue of Theorem \ref{main2},
				we have  a long exact sequence in cohomology:{\footnotesize
					$$
					\dots\rightarrow \cohom^p(\Omega_c(G/H)_\Gamma)\overset{\cohom(\iota)}{\longrightarrow} \cohom_c^p(G/H)\overset{\cohom(m)}{\longrightarrow} \cohom^p(\Omega(G/H)^\Gamma)\overset{\delta}{\longrightarrow}
					\cohom^{p+1}(\Omega_c(G/H)_\Gamma)\rightarrow\dots
					$$}
				In the particular case where $H$ is a maximal compact subgroup of a connected Lie group $G$, the manifold $G/H$ is then a contractible, which leads to:
				$$\cohom^p(\Omega_c(G/H)_\Gamma)\simeq \cohom^{p-1}(\Omega(G/H)^\Gamma)$$ and when the action of $\Gamma$ on $G/H$ is free, we have $\cohom^p(\Omega_c(G/H)_\Gamma)\simeq \cohom^{p-1}(\Gamma\backslash G/H).$ \end{enumerate}
			
		\end{Example}

		%<-------------------
		
		% Your bilbigraphy           %<-------------------

			A. Abouqateb\\               %<-------------------
			Cadi Ayyad University, Faculty of Sciences and Technologies, Department of Mathematics B.P.549 Gueliz Marrkesh. Morocco\\            %<-------------------
			a.abouqateb@uca.ac.ma                %<-------------------

			M. Boucetta\\               %<-------------------
			Cadi Ayyad University, Faculty of Sciences and Technologies, Department of Mathematics B.P.549 Gueliz Marrkesh. Morocco\\            %<-------------------
			m.boucetta@uca.ac.ma                %<-------------------

			M. Nabil\\               %<-------------------
			Cadi Ayyad University, Faculty of Sciences and Technologies, Department of Mathematics B.P.549 Gueliz Marrkesh. Morocco\\            %<-------------------
			mehdi1nabil@gmail.com                %<-------------------

\end{document}